\documentclass[12pt, reqno]{amsart}
\usepackage{amsmath, amsthm, amscd, amsfonts, amssymb, graphicx, color}
\usepackage[bookmarksnumbered, colorlinks, plainpages]{hyperref}
\hypersetup{colorlinks=true,linkcolor=red, anchorcolor=green, citecolor=cyan, urlcolor=red, filecolor=magenta, pdftoolbar=true}

\newtheorem{theorem}{Theorem}[section]
\newtheorem{lemma}[theorem]{Lemma}
\newtheorem{proposition}[theorem]{Proposition}
\newtheorem{corollary}[theorem]{Corollary}
\theoremstyle{definition}
\newtheorem{definition}[theorem]{Definition}
\newtheorem{example}[theorem]{Example}

\theoremstyle{remark}
\newtheorem{remark}[theorem]{Remark}
\numberwithin{equation}{section}

\begin{document}

\setcounter{page}{1}

\title[BSE-property]{BSE-property for some certain Segal and Banach algebras}

\author[M. Fozouni \MakeLowercase{and} M. Nemati ]{Mohammad Fozouni$^1$$^{*}$ \MakeLowercase{and} Mehdi Nemati$^2$}

\address{$^{1}$Department of Mathematics, Gonbad Kavous University, P.O.Box 163, Gonbad-e Kavous, Iran.}
\email{\textcolor[rgb]{0.00,0.00,0.84}{fozouni@gonbad.ac.ir}}

\address{$^{2}$Department of Mathematical Sciences, Isfahan University of Technology, Isfahan 84156-83111, Iran;
\newline
}
\email{\textcolor[rgb]{0.00,0.00,0.84}{m.nemati@cc.iut.ac.ir}}





\date{Received: xxxxxx; Revised: yyyyyy; Accepted: zzzzzz.
\newline \indent $^{*}$Corresponding author}

\begin{abstract}
For a commutative semi-simple Banach algebra ${A}$ which is an ideal in its second dual  we  give a necessary and sufficient condition for
an essential abstract Segal algebra in ${A}$ to be a BSE-algebra. We show that
 a large class of abstract Segal algebras in the Fourier algebra $A(G)$ of a locally compact group $G$ are BSE-algebra if and only if they have bounded weak  approximate identities. Also, in the case that $G$ is discrete we show that $A_{\rm cb}(G)$ is a BSE-algebra if and only if $G$ is weakly amenable.
 We study the BSE-property of some certain Segal algebras implemented by local functions that were
recently introduced by J. Inoue and S.-E. Takahasi. Finally we give a similar construction for the group algebra implemented by a measurable and sub-multiplicative function.
\end{abstract}

\maketitle
\section{Introduction and Preliminaries}

Let $G$ be a locally compact abelian group. A subspace $\mathcal{S}$ of $L^{1}(G)$ is called a (Reiter) Segal algebra if it satisfies the following conditions:
\begin{enumerate}
  \item $\mathcal{S}$ is dense in $L^{1}(G)$.
  \item $\mathcal{S}$ is a Banach space under some norm $\|\cdot\|_{\mathcal{S}}$  such that $\|f\|_{1}\leq \|f\|_{\mathcal{S}}$ for each $f\in \mathcal{S}$.
  \item $L_{y}f$ is in $\mathcal{S}$ and $\|f\|_{\mathcal{S}}=\|L_{y}f\|_{\mathcal{S}}$ for all $f\in \mathcal{S}$ and $y\in G$ where $L_{y}f(x)=f(y^{-1}x)$.
    \item For all $f\in \mathcal{S}$, the mapping $y\longrightarrow L_{y}f$ is continuous with respect to $\|\cdot\|_{\mathcal{S}}$.
\end{enumerate}
  In \cite{Burnham}, J. T. Burnham with changing $L^{1}(G)$ by  an arbitrary Banach algebra $A$, gave a generalization of Segal algebras and introduced the notion of an abstract Segal algebra.
Recall that a Banach
algebra $B$ is an abstract Segal algebra of  a Banach
algebra $A$ if

\begin{enumerate}
  \item $B$ is a dense left ideal in $A$,
  \item there exists $M > 0$ such that $\|b\|_{A}\leq M \|b\|_{B}$ for each $b\in B$,
  \item there exists $C > 0$ such that $\|ab\|_{B}\leq
C\|a\|_{A} \|b\|_{B}$ for each $a, b\in B$.
\end{enumerate}

Recently, J. Inoue and S.-E. Takahasi in \cite{Inoue} investigated abstract Segal algebras in a non-unital commutative semi-simple regular Banach algebra ${A}$ such that ${A}$ has a bounded approximate identity in ${A}_{c}$ where
\begin{equation*}
{A}_{c}=\{a\in {A} : \widehat{a} \text{ has compact support}\},
\end{equation*}
and $\widehat{a}$ denotes the Gel'fand transform of $a$. It is well-known that $L^{1}(G)$ is a commutative semi-simple regular Banach algebra with a bounded approximate identity with compact support; see \cite{Kaniuth2}.

A commutative Banach algebra $A$ is  without order if  for $a\in A$, the condition $aA=\{0\}$ implies $a=0$ or equivalently $A$ does not have any non-zero annihilator. For example if ${A}$ has an approximate identity, then it is without order.
A linear operator $T$ on $A$ is called a multiplier if it satisfies $aT(b)=T(a)b$ for all $a, b\in {A}$. Suppose that $\mathcal{M}(A)$ denotes the space of all multipliers of the Banach algebra $A$.
If $\Delta(A)$ denotes the space of all characters of $A$; that is, non-zero homomorphisms from
$A$ into $\mathbb{C}$, then for each $T\in \mathcal{M}(A)$, there exists a unique bounded continuous function $\widehat{T}$ on $\Delta(A)$ such that $\widehat{T(a)}(\phi)=\widehat{T}(\phi)\widehat{a}(\phi)$ for all $a\in A$ and $\phi\in \Delta(A)$; see \cite[Proposition 2.2.16]{Kaniuth2}. Let $\widehat{\mathcal{M}(A)}$ denote the space of all $\widehat{T}$ corresponding to $T\in \mathcal{M}(A)$.

A bounded continuous function $\sigma$ on $\Delta(A)$ is called a BSE-function if there exists a constant $C>0$ such that for each $\phi_{1},...,\phi_{n}\in\Delta(A)$ and complex numbers $c_{1},...,c_{n}$, the inequality
\begin{equation*}
\left|\sum_{i=1}^{n}c_{i}\sigma(\phi_{i})\right|\leq C\left\|\sum_{i=1}^{n}c_{i}\phi_{i}\right\|_{A^{*}}
\end{equation*}
holds. Let $C_{\mathrm{BSE}}(\Delta(A))$ be the set of all $\mathrm{BSE}$-functions.

A without order commutative Banach algebra $A$ is called a BSE-algebra if $$C_{\mathrm{BSE}}(\Delta(A))=\widehat{\mathcal{M}(A)}.$$

The theory of  BSE-algebras was first introduced and investigated by Takahasi and Hatori; see \cite{Takahashi} and two other notable works \cite{Kaniuth, Kamali2}. In \cite{Kamali2}, the authors answered a question raised in \cite{Takahashi}. Examples of BSE-algebras are the group algebra $L^1(G)$ of a locally compact abelian
group $G$, the Fourier algebra $A(G)$ of a locally compact amenable
group $G$, all commutative $C^*$-algebras, the disk algebra, and the Hardy algebra on the
open unit disk.

A net $\{a_{\alpha}\}$ in $A$ is called a bounded weak approximate identity (BWAI) for $A$ if $\{a_{\alpha}\}$ is bounded in $A$ and
\begin{equation*}
\lim_{\alpha}\phi(a_{\alpha}a)=\phi(a)\hspace{0.5cm}(\phi\in \Delta(A), a\in A),
\end{equation*}
or equivalently, $\lim_{\alpha}\phi(a_{\alpha})=1$ for each $\phi\in \Delta(A)$. Clearly, each BAI of $A$ is a BWAI and the converse is not valid in general; see \cite{Jones} and \cite{Laali}. Note that  bounded weak approximate identities are  important to decide whether a commutative Banach algebra is a BSE-algebra or not. For example, if $A$ is a Banach algebra then $\widehat{\mathcal{M}(A)}\subseteq C_{\mathrm{BSE}}(\Delta(A))$ if and only if $A$ has a BWAI; \cite[Corollary 5]{Takahashi}. Also, it was shown in \cite{ino2} that
a Segal algebra $S(G)$ on a locally compact abelian group $G$
is a BSE-algebra if and only if $S(G)$ has a  BWAI.

For undefined concepts and notations appearing in the sequel, one can consult \cite{Dales, Kaniuth2}.

The outline of the next sections is as follows:

In \S 2, for a commutative semi-simple Banach algebra ${A}$ which is an ideal in its second dual  we  give a necessary and sufficient condition for
an essential abstract Segal algebra
in ${A}$ to be a BSE-algebra. We show that
a large class of abstract Segal algebras in the Fourier algebra $A(G)$ of a locally compact group $G$ are BSE-algebras if and only if they have  bounded weak approximate identities. Also, for discrete groups $G$, we show that $A_{\rm cb}(G)$ is a BSE-algebra if and only if $G$ is weakly amenable.

 In \S 3, we study the BSE-property of the Segal algebra $A_{\tau(n)}$ in $A$ which was introduced by Inoue and Takahasi and in the case that $(A, \|\cdot\|_{X})$ is a BSE-algebra, we show that $A_{\tau(n)}$ is a BSE-algebra if and only if $\tau$ is bounded where $\tau: X\longrightarrow \mathbb{C}$ is a certain continuous function.
Also, we compare the BSE-property between
$A$ and $A_{\tau(n)}$.
 In \S 4, motivated by the definition of $A_{\tau(n)}$, for an arbitrary locally compact (abelian) group $G$, and a measurable sub-multiplicative function $\tau: G\longrightarrow \mathbb{C}^{\times}$, we define the Banach algebra $L^{1}(G)_{\tau(n)}$. Then we  investigate the BSE-property of this algebra.
\section{BSE-Abstract Segal Algebras}
Suppose that $B$ is an abstract Segal algebra of Banach algebra $A$. Endow $\Delta({A})$ and $\Delta({B})$ with the
Gel'fand topology,
the map $\varphi\mapsto\varphi|_{B}$ is a homeomorphism from
$\Delta({A})$ onto $\Delta({B})$; see \cite[Theorem 2.1]{Burnham}.

Let ${B}$ be an abstract Segal algebra with respect to ${A}$.
We say that ${B}$ is essential if $\langle {A}{B}\rangle $ is $\|\cdot\|_{B}$--dense in ${B}$, where $\langle {A}{B}\rangle$ denotes the linear span of  ${A}{B}=\{ab: a\in {A}, b\in {B}\}$.
\begin{theorem}\label{seg1}
Let ${A}$  be a semi-simple commutative Banach algebra which is an
ideal in its second dual ${A}^{**}$. Suppose that
${B}$ is an essential abstract Segal algebra
in ${A}$. Then the following statements are equivalent.

\rm{(i)} ${B}$ is a BSE-algebra.

\rm{(ii)} ${B}={A}$ and ${A}$ is a BSE-algebra.
\end{theorem}
\begin{proof}

Suppose that  ${B}$ is a BSE-algebra. Then  by \cite[Corollary 5]{Takahashi}, ${B}$ has a BWAI, say $(g_\gamma)_\gamma$. It is clear that $(g_\gamma)_\gamma$ is also a BWAI for ${A}$. So, by \cite[Theorem 3.1]{Kaniuth} ${A}$ is a BSE-algebra and has a bounded approximate identity, say $(e_\alpha)_\alpha$. By the Cohen factorization theorem $AB$ is a closed linear subspace of $B$. But $B$ is essential, so $B=AB$. Now, let $b\in{B}$. Then   $b = ca$ for some $c\in{B}$ and $a\in{A}$. Since ${A}$ is an ideal in its second dual, it follows that the operator ${\rho}_a: {A} \rightarrow {A}$ defined by $\rho_a(a')=aa',$  $(a'\in {A})$ is weakly compact; see \cite[Lemma 3]{dh}. Letting $\iota : B\hookrightarrow A$ and $L_c : A\rightarrow B : a^{'}\rightarrow ca^{'}$, $\rho_b =L_c\circ \rho_a\circ \iota$ is weakly compact by \cite[Proposition VI.5.2]{Conway}. Therefore, ${B}$ is an ideal in its second dual. Since by \cite[Theorem 2.1]{Burnham} ${B}$ is semi-simple,  \cite[Theorem 3.1]{Kaniuth} implies that ${B}$ has a bounded approximate identity. Thus ${A}={B}$ by \cite[Theorem 1.2]{Burnham}. That (ii) implies (i) is trivial.
\end{proof}
\begin{example}\label{seg2}
Let $G$ be a locally compact group and let  $A(G)$ be the Fourier algebra  of $G$.
It was shown in \cite[Theorem 5.1]{Kaniuth} that $A(G)$ is a BSE-algebra if and only if $G$ is amenable. Moreover,   $A(G)$ is an ideal in its second dual if and only if  $G$ is discrete; see \cite[Lemma 3.3]{F91}. Therefore, by Theorem \ref{seg1} if $G$ is discrete, then an
 essential abstract Segal algebra $SA(G)$
in $A(G)$ is a BSE-algebra if and only if  $SA(G)=A(G)$ and $G$ is amenable.
\end{example}
Let $G$ be a locally compact group and let  $L^r(G)$ be the Lebesgue $L^r$-space  of $G$, where $1\leq r<\infty$. Then
$$
SA^r(G):= L^r(G)\cap A(G)
$$
with the norm
$||| f ||| =\|f \|_r + \|f \|_{A(G)}$ and the pointwise
product is an abstract Segal algebra in $A(G)$.

\begin{corollary}
Let $G$ be a discrete group and let $1\leq r\leq 2$. Then $SA^r(G)$ is a \rm{BSE}-algebra if and only if $G$ is finite
\end{corollary}
\begin{proof}
If $G$ is finite, then $SA^r(G)=A(G)$. So, the result follows from Theorem \ref{seg1}.

For the converse, first note that $SA^r(G)=l^r(G)$ and the norms $\|\cdot\|_r$ and $|||\cdot |||$ on $SA^r(G)$ are
equivalent by the open mapping theorem. In fact, $l^2(G)\subseteq \delta_e*l^2(G)\subseteq A(G)$, where $\delta_e$ is the point mass at the identity element $e$ of $G$. So,  if  $1\leq r\leq 2$, then $l^r(G)\subseteq l^2(G)$ and $$l^r(G)= l^r(G)\cap l^2(G)\subseteq SA^r(G)\subseteq l^r(G).$$ Moreover, it is clear that  $l^r(G)$  has an approximate identity and consequently it
is an essential abstract Segal algebra in $A(G)$. Therefore, if
$SA^r(G)=l^r(G)$ is a BSE-algebra, then $A(G)=l^r(G)$ by Example \ref{seg2}. Thus, $A(G)=l^2(G)$, is a reflexive predual of a $W^*$-algebra. This
implies, as is  known, that $A(G)$ is finite dimensional; see \cite{Sakai}. Thus $G$ is finite, which completes  the proof.
\end{proof}
For a locally compact group $G$, we recall that $A(G)$ is always an ideal in the  Fourier-Stieltjes algebra $B(G)$ and note that $\mathcal{M}(A(G))=B(G)$ when $G$ is amenable.  The spectrum of $A(G)$ can be canonically
 identified with $G$. More precisely, the map $x\longrightarrow
 \varphi_x$  where $\varphi_x(u) = u(x)$ for all $u\in A(G)$ is a
 homeomorphism from $G$ onto $\Delta(A(G))$. Note that by \cite[Theorem 2.1]{Burnham}, $\Delta(SA(G))= \Delta(A(G))=G$.
\begin{theorem}\label{seg3}
Let $G$  be a locally compact group and let
$SA(G)$ be an  abstract Segal algebra
in $A(G)$ such that $B(G)\subseteq \mathcal{M}(SA(G))$, i.e., for each $u\in B(G),$  we have $uSA(G)\subseteq SA(G)$. Then $SA(G)$ is a \rm{BSE}-algebra if and only if
 $SA(G)$ has a BWAI.
\end{theorem}
\begin{proof}
By \cite[Corollay 5]{Takahashi}, if  $SA(G)$ is a BSE-algebra, then  it has a BWAI.

Conversely,  suppose that $SA(G)$ has a BWAI, say $(e_\gamma)_\gamma$. Then $\widehat{\mathcal{M}(SA(G))}\subseteq C_{BSE}(\Delta(SA(G)))$  by  \cite[Corollary 5]{Takahashi}. Moreover, it is clear that $(e_\gamma)_\gamma$  is also a BWAI for $A(G)$. Consequently,  we conclude that  $G$ is amenable by  \cite[Theorem 5.1]{Kaniuth}.
Now, we need to show  the reverse inclusion. Since $SA(G)$  is an abstract Segal algebra in $A(G)$, there exists $M>0$ such that $\|u\|_{A(G)}\leq M\|u\|_{SA(G)}$ for all $u\in SA(G)$. Thus, for any $x_1,...,x_n\in G$ and $c_1,...,c_n\in{\Bbb C}$,
$$
\left\|\sum_{j=1}^n c_j\varphi_{x_j}\right\|_{SA(G)^*}\leq
M \left\|\sum_{j=1}^n c_j\varphi_{x_j}\right\|_{A(G)^*}.
$$
This implies that
\begin{eqnarray*}
C_{BSE}(\Delta(SA(G)))&\subseteq& C_{BSE}(\Delta(A(G)))\\
&=&\widehat{B(G)}\\
&\subseteq&
\widehat{\mathcal{M}(SA(G))}.
\end{eqnarray*}
Hence, $SA(G)$ is a BSE-algebra.
\end{proof}
\begin{example}
(1) Let $G$ be a locally compact group and let $1\leq r<\infty$. Now, since $\|u\|_\infty \leq\|u\|_{B(G)}$ for all $u\in B(G)$, it follows that $uL^r(G)\subseteq L^r(G)$. This implies that  $B(G)\subseteq \mathcal{M}(SA^r(G))$. Thus $SA^r(G)$
is a BSE-algebra if and only if it has a BWAI .

(2) Let $S_0(G)$ be
the Feichtinger Segal algebra in $A(G)$. Then $B(G)\subseteq \mathcal{M}(S_0(G))$; see \cite[ Corollary 5.2]{samei}. Thus $S_0(G)$
is a BSE-algebra if and only if it has a BWAI
\end{example}
\begin{corollary}\label{seg4}
Let $G$  be a locally compact group and let
$SA(G)$ be an essential  abstract Segal algebra
in $A(G)$. Then $SA(G)$ is a \rm{BSE}-algebra if and only if
 $SA(G)$ has a BWAI.
\end{corollary}
\begin{proof}
Suppose that $SA(G)$ has a BWAI. Hence $G$ is amenable, therefore $A(G)$ has a bounded approximate identity. Now, by the same argument as the proof of  Theorem \ref{seg1} we can show that
  $SA(G)=A(G)SA(G)$.  Consequently, for each $u\in B(G)$
    $$uSA(G)=uA(G)SA(G)\subseteq A(G)SA(G)=SA(G),$$ which implies that $B(G)\subseteq \mathcal{M}(SA(G))$. Hence, $SA(G)$ is a BSE-algebra by Theorem \ref{seg3}.
\end{proof}
Suppose that $G$ is a locally compact group and $\mathcal{M}_{\rm cb}A(G)$ denotes the Banach algebra of completely bounded multipliers of $A(G)$, that is, the continuous and bounded functions $\nu$ on $G$ such that $\nu A(G)\subseteq A(G)$ and the map $L_{\nu}$ defined by $L_{\nu}(u)=\nu u$ is completely bounded; see \cite{ER} for a complete course on operator space theory. Note that, by $A(G)=VN(G)_*$, $\mathcal{M}_{\rm cb}A(G)$ is a completely contractive Banach algebra, where $VN(G)$ is the group von Neumann algebra. It is well-known that $A(G)\subseteq B(G)\subseteq \mathcal{M}_{\rm cb}A(G)$.
 Now, let
$$A_{\rm cb}(G)=\overline{A(G)}^{\|\cdot\|_{\mathcal{M}_{\rm cb}A(G)}}.$$

This algebra was first introduced by Forrest in \cite{For2}. See also some recent works \cite{Fore, FRS, Nem} for more details and properties.

We end this section with the following result regarding the BSE-property of $A_{\rm cb}(G)$ in the case that $G$ is discrete. Recall that a locally compact group $G$ is said to be weakly amenable if $A(G)$ has an approximate identity which is bounded in $\|\cdot\|_{A_{\rm cb}(G)}$ or equivalently  $A_{\rm cb}(G)$ has a bounded approximate identity; see \cite[Proposition 1]{Fore}.

\begin{theorem}\label{Th3}
  Suppose that $G$ is a discrete locally compact group. Then $A_{\rm cb}(G)$ is a BSE-algebra if and only if $G$ is weakly amenable.
   \end{theorem}
   \begin{proof}
By \cite[Lemma 4.1]{Nem} we know that $G$ is discrete if and only if $A_{\rm cb}(G)$ is a closed ideal in its second dual. Also, it is clear that  $A_{\rm cb}(G)$ is commutative and semi-simple. Now the result follows by \cite[Theorem 3.1]{Kaniuth}.
   \end{proof}
\section{Segal Algebras Implemented by Local Functions}
In  this section we focus on a certain Segal algebra that was recently introduced by Inoue and Takahasi. Let $X$ be a non-empty locally compact Hausdorff space. A subalgebra $A$ of $C_{0}(X)$ is called a Banach function algebra if $A$ separates strongly the points of $X$ (that is, for each $x, y\in X$ with $x\neq y$, there exists $f\in A$ such that $f(x)\neq f(y)$ and for each $x\in X$, there exists $f\in A$ such that $f(x)\neq 0$) and with a norm $\|\cdot\|$, $(A, \|\cdot\|)$ is a Banach algebra.

Suppose that $(A, \|\cdot\|)$ is a natural  regular Banach function algebra on a locally compact, non-compact Hausdorff space $X$ with a bounded approximate identity $\{e_{\alpha}\}$ in $A_{c}$.  We  recalling the following  definitions from \cite{Inoue}.
\begin{definition}
A complex-valued continuous function $\sigma$ on $X$ is called a local $A$-function if for all $f\in A_{c}$, $f\sigma\in A$. The set of all local $A$-functions is denoted by $A_{\mathrm{loc}}$
\end{definition}
\begin{definition}
For positive integer $n$ and a continuous complex-valued function $\tau$ on $X$, put
\begin{align*}
A_{\tau(n)}&= \left\{f\in A : f\tau^{k}\in A\quad(0\leq k\leq n)\right\},\\
\|f\|_{\tau(n)}&=\sum_{k=0}^{n}\|f\tau^{k}\|.
\end{align*}
\end{definition}
 In the sequel of this section, suppose that $n$ is a constant positive integer and $\tau\in A_{\mathrm{loc}}$.

By \cite[Theorem 5.4]{Inoue}, if $\tau\in A_{\mathrm{loc}}$, then $(A_{\tau(n)}, \|\cdot\|_{\tau(n)})$ is a Segal algebra in $A$ such that $\Delta(A_{\tau(n)})=\Delta(A)=X$, that is, $x\longrightarrow \phi_{x}$ is a homeomorphism from $X$ onto $\Delta(A_{\tau(n)})$.

Also, one can see  that $A_{\tau(n)}$ is a Banach function algebra, because for each $x, y\in X$ with $x\neq y$, there exists $f\in A_{\tau(n)}$ such that  $f(x)=\phi_{x}(f)\neq \phi_{y}(f)=f(y)$ and by using the Urysohn lemma for each $x\in X$, there exists $f\in A_{\tau(n)}$ with $f(x)\neq 0$. Note that by \cite[Theorem 3.5]{Inoue}, $A_{\tau(n)}$ is  Tauberian. Recall that a Banach algebra $A$ is Tauberian if $A_{c}$ is dense in $A$.


The following theorem is one of our main results in this section.
\begin{theorem}\label{Th} Suppose that $(A, \|\cdot\|_{X})$ is a $\mathrm{BSE}$-algebra where $\|\cdot\|_{X}$ is the uniform norm. Then the following statements are equivalent.

{\rm (i)} $A_{\tau(n)}$ is a $\mathrm{BSE}$-algebra.

   {\rm (ii)} $\tau$ is bounded.

\end{theorem}
\begin{proof}
(i) $\rightarrow$(ii).
Suppose that   $A_{\tau(n)}$  is a BSE-algebra. Then it has a BWAI, say $\{f_\alpha\}$. So, there exists a constant  $M>0$ such that
$$\|f_{\alpha}\|_{{\tau(n)}}<M,\quad \lim_{\alpha}f_{\alpha}(x)=1\quad\quad(x\in X).$$
On the other hand, $\|f_{\alpha}\tau\|_{X}\leq \|f_{\alpha}\|_{{\tau(n)}}$, hence we have
$$|\tau(x)|\leq M\hspace{0.5cm}(x\in X).$$
Therefore, $\tau$ is bounded.

(ii)$\rightarrow $(i).
Let $\tau$ be bounded by $M$, that is, $\|\tau\|_{X}<M$. Clearly, $A_{c}\subseteq A_{\tau(n)}$. For each $f\in A$, there exists a net $\{f_{\alpha}\}$ in $A_{c}$ such that $\|f_{\alpha}-f\|_{X}\longrightarrow 0$. Now, $\{f_{\alpha}\tau\}$  is in $A$ and $f\tau=\lim_{\alpha}f_{\alpha}\tau$, since
\begin{equation*}
\|f\tau-f_{\alpha}\tau\|_{X}=\|(f-f_{\alpha})\tau\|_{X}\leq M\|f-f_{\alpha}\|_{X}\longrightarrow 0.
\end{equation*}
 So, $f\tau\in A$.   Similarly, one can see that $f\tau^{k}$ for each $1<k\leq n$ is in $A$. Therefore, $f\in A_{\tau(n)}$. Hence $A=A_{\tau(n)}$. Finally, for each $f\in A$
 $$\|f\|_{X}\leq \|f\|_{\tau(n)}\leq \|f\|_{X}\left(\sum_{k=0}^{n}M^{k}\right).$$
  Therefore, $A_{\tau(n)}$ is topologically isomorphic to $A$, so $A_{\tau(n)}$ is a BSE-algebra.
\end{proof}
Recall that $\tau\in A_{\mathrm{loc}}$ is called a rank $\infty$ local $A$-function, if for each $k=0,1,2,\cdots$, the inclusion $A_{\tau(k)}\supsetneqq A_{\tau(k+1)}$ holds. By \cite[Proposition 8.2 (ii)]{Inoue}, if $\|\tau\|_{X}=\infty$, then $\tau$ is a rank $\infty$ local $A$-function.

As an application of the above theorem, we give the following result which provide for us examples of Banach algebras without any BWAI.
\begin{corollary}
Let $A=C_{0}(\mathbb{R})$  and $\tau(x)=x$ for every $x\in \mathbb{R}$. Then
\begin{equation}\label{eq}
A\supsetneqq A_{\tau(1)}\supsetneqq A_{\tau(2)}\supsetneqq \cdots \supsetneqq A_{\tau(n)}\supsetneqq \cdots.
\end{equation}
For each $k=1, 2, 3, \cdots$, $A_{\tau(k)}$ is not a $\mathrm{BSE}$-algebra and  has no BWAI .
\end{corollary}
\begin{proof}
By \cite[Theorem 3]{Takahashi}, if $A=C_{0}(\mathbb{R})$, we know that $A$ is a BSE-algebra. Since for each $x\in \mathbb{R}$, there exists $f\in A$ such that $f=\tau$ on a neighborhood of $x$, hence by \cite[Proposition 7.2]{Inoue}, $\tau$ is an element of $A_{\mathrm{loc}}$. But $\tau$ is not bounded and this implies that $A_{\tau(k)}$ is not a BSE-algebra by  Theorem \ref{Th}. Also, since  $\tau$ is not bounded, $\tau$ is a rank $\infty$ local $A$-function, hence we have equation \ref{eq}.
 Finally, if $A_{\tau(k)}$ has a BWAI, then similar to the proof of Theorem \ref{Th} one has $\|\tau\|_{\mathbb{R}}< \infty$ which is impossible.
\end{proof}
By the above corollary if $A$ is a BSE-algebra, then $A_{\tau(n)}$ is not necessarily a BSE-algebra. For the converse we have the following proposition.
\begin{proposition}\label{Th1}
Suppose that $\tau\in A_{\mathrm{loc}}$  and $A$ is an ideal in its second dual.
If $A_{\tau(n)}$ is a $\mathrm{BSE}$-algebra, then $A$ is a $\mathrm{BSE}$-algebra.
\end{proposition}
\begin{proof}
If $\{e_{\alpha}\}$ is a b.a.i. for $A$, then it follows from our definition of $A_{\tau(n)}$ that $\{e_{\alpha}\}$ is also an approximate identity for $A_{\tau(n)}$ when viewed as a Banach $A$-module.

So, $A_{\tau(n)}$ is an essential abstract Segal algebra with respect to $A$, because for each $b\in A_{\tau(n)}$, $\|b-e_\alpha b\|_{\tau(n)}\rightarrow 0$. Now, using Theorem \ref{seg1} we have the result.
\end{proof}
We do not know whether Proposition \ref{Th1} fails  if the assumption that $A$ is an ideal in its second dual is dropped.
\begin{remark} If $X$ is a discrete space, then $A$ is an ideal in its second dual. Since $A$ is a semi-simple commutative and Tauberian Banach algebra, using \cite[Remark 3.5]{Kaniuth}, we conclude that $A$ is an ideal in its second dual.
\end{remark}

\section{A Construction on Group Algebras}
Let $G$ be a locally compact group and let
$L^{1}(G)$ be the space of all measurable and integrable complex-valued functions (equivalent classes with respect to the almost everywhere equality relation) on $G$ with respect to the left Haar measure of $G$. The convolution product  of the functions $f$ and $g$ in $L^1(G)$ is defined by
\begin{equation*}
f\ast g(x)=\int_{G}f(y)g(y^{-1}x)dy.
\end{equation*}
For each $f\in L^{1}(G)$, let $\|f\|_{1}=\int_{G}|f(x)|dx$.
It is well-known that $L^{1}(G)$ endowed with the norm $\|\cdot\|_{1}$ and the convolution product is a Banach algebra called the group algebra of $G$; see \cite[Section 3.3]{Dales} for more details.

The following lemma was proved by Bochner and Schoenberg for $G=\mathbb{R}$ in (1934) and by Eberlein for general locally compact abelian (LCA) groups in (1955). Here, we explicitly observe that a result due to E. Kaniuth and
A. \"{U}lger gives another proof.

Recall that for a LCA group $G$, the dual group of $G$, $\widehat{G}$ is defined as the set of all continuous homomorphisms from $G$ to $\mathbb{T}$ where $\mathbb{T}=\{z\in \mathbb{C} : |z|=1\}$. It is well-known that $\widehat{G}$ is a LCA group with the pointwise-defined operation.

\begin{lemma}\label{BSE}
Suppose that $G$ is $\mathrm{LCA}$ group. Then $L^{1}(G)$ is a $\mathrm{BSE}$-algebra.
\end{lemma}
\begin{proof}
We know that $L^{1}(G)$ is isometrically isomorphic to $A(\widehat{G})$. But   $\widehat{G}$ is amenable. Therefore, by \cite[Theorem 5.1]{Kaniuth} $A(\widehat{G})$, and so $L^{1}(G)$, is a BSE-algebra.
\end{proof}
In the sequel, motivated by the construction of $A_{\tau(n)}$ in the preceding section, we introduce a subalgebra of the group algebra $A=L^{1}(G)$ where $G$ is a locally compact group.

Recall that $\varphi: G\longrightarrow \mathbb{C}^{\times}$ is sub-multiplicative if
$$|\varphi(xy)|\leq |\varphi(x)\|\varphi(y)|\quad (x, y\in G),$$
where $\mathbb{C}^{\times}$ denotes the multiplicative group of non-zero complex numbers.

For a measurable sub-multiplicative function $\tau: G\longrightarrow \mathbb{C}^{\times}$ and each $n\in {\Bbb N}$, put
\begin{align*}
L^{1}(G)_{\tau(n)}&=\{f\in L^{1}(G) : f\tau,\ldots,f\tau^{n}\in L^{1}(G)\}\\
\|f\|_{\tau(n)}&=\sum_{k=0}^{n}\|f\tau^{k}\|_{1}\hspace{0.5cm}(f\in L^{1}(G)_{\tau(n)}).
\end{align*}
As the first result in this section, we show that $L^{1}(G)_{\tau(n)}$ is a Banach algebra as follows.
\begin{proposition}
$L^{1}(G)_{\tau(n)}$ is a Banach algebra with the convolution product and the norm $\|\cdot\|_{\tau(n)}$.
\end{proposition}
\begin{proof}
For each $f, g\in L^{1}(G)_{\tau(n)}$ and $1\leq k\leq n$, we have
\begin{align*}
\|(f\ast g)\tau^{k}\|_{1}&\leq \int\int |f(y)\|g(y^{-1}x)\|\tau^{k}(x)|dydx\\
&=\int\int |f(y)\|g(y^{-1}x)\|\tau^{k}(x)|dxdy\\
&=\int\int |f(y)\|g(x)\|\tau^{k}(yx)|dxdy\\
&\leq\|f\tau^{k}\|_{1}\|g\tau^{k}\|_{1}.
\end{align*}
For each $1\leq k\leq n$, since $\tau$ is measurable, $(f\ast g)\tau^{k}$ is measurable and by the above inequality  $\|(f\ast g)\tau^{k}\|_{1}<\infty$. Therefore, $f\ast g$ is in $L^{1}(G)_{\tau(n)}$. Also, we have
\begin{align*}
\|f\ast g\|_{\tau(n)}&=\|f\ast g\|_{1}+\sum_{k=1}^{n}\|(f\ast g)\tau^{k}\|_{1}\\
&\leq \|f\|_{1}\|g\|_{1}+\sum_{k=1}^{n}\|f\tau^{k}\|_{1}\|g\tau^{k}\|_{1}\\
&\leq \|f\|_{\tau(n)}\|g\|_{\tau(n)}.
\end{align*}
To see the completeness of $\|\cdot\|_{\tau(n)}$, let $\{f_{i}\}$ be a Cauchy sequence in $L^{1}(G)_{\tau(n)}$. So, there exist $f\in L^{1}(G)$ and $g_{k}\in L^{1}(G)$ for each $1\leq k\leq n$ such that
\begin{equation*}
\lim_{i\rightarrow \infty}\|f_{i}-f\|_{1}=0,\quad \lim_{i\rightarrow \infty}\|f_{i}\tau^{k}-g_{k}\|_{1}=0.
\end{equation*}
Since $\lim_{i\rightarrow \infty}\|f_{i}-f\|_{1}=0$, there exists a subsequence $\{f_{i_{m}}\}$ such that $$\lim_{i_{m}}f_{i_{m}}(x)=f(x) \text{ a.e., for } x\in G.$$ Also, since $\lim_{i_{m}}\|f_{i_{m}}\tau^{k}-g_{k}\|_{1}=0$, there exists a subsequence $\{f_{i_{m , k}}\}$ of $\{f_{i_{m}}\}$ such that $$\lim_{i_{m,k}}f_{i_{m,k}}(x)\tau^{k}(x)=g_{k}(x) \text{ a.e., for } x\in G.$$
Therefore, $f\tau^{k}=g_{k}$ a.e., so, $f$ is an element of $L^{1}(G)_{\tau(n)}$ such that
\begin{align*}
\|f_{i}-f\|_{\tau(n)}&=\|f_{i}-f\|_{1}+\|f_{i}\tau-f\tau\|_{1}+\cdots+\|f_{i}\tau^{n}-f\tau^{n}\|_{1}\\
&=\|f_{i}-f\|_{1}+\|f_{i}\tau-g_{1}\|_{1}+\cdots+\|f_{i}\tau^{n}-g_{n}\|_{1}\longrightarrow 0.
\end{align*}
Hence, $(L^{1}(G)_{\tau(n)},\|\cdot\|_{\tau(n)})$ is complete.
\end{proof}
In the sequel, we suppose that $G$ is a LCA group.
\begin{theorem}\label{Th2}
If $\tau$ is bounded, then $L^{1}(G)_{\tau(n)}$ is a $\mathrm{BSE}$-algebra.
\end{theorem}
\begin{proof}
If $\tau$ is bounded by $M$, then for each $f\in L^{1}(G)$ and $1\leq k\leq n$ we have $f\tau^{k}\in L^{1}(G)$ and
$$\|f\|_{1}\leq \|f\|_{\tau(n)}=\sum_{k=0}^{n}\|f\tau^{k}\|_{1}\leq \|f\|_{1}(\sum_{k=0}^{n}M^{k}).$$

So, $L^{1}(G)_{\tau(n)}$ is topologically isomorphic to $L^{1}(G)$ and hence it is a BSE-algebra by Lemma \ref{BSE} and \cite[Corollary 1.3]{Kaniuth}.
\end{proof}
When $\tau$ satisfying $|\tau(x)|\geq 1$ a.e., for $x\in G$, we show in the sequel that $L^{1}(G)_{\tau(n)}$ is in fact a Beurling algebra, the definition of which we now recall.

A weight on $G$ is a measurable  function $w: G\longrightarrow (0,\infty)$ such that  $w(xy)\leq w(x)w(y)$ for all $x, y\in G$. The Beurling algebra $L^{1}(G, w)$ is defined to be the space of all measurable complex-valued functions $f$ on $G$ such that $\|f\|_{1,w}=\int|f(x)|w(x)\ dx<\infty$. The Beurling algebra with the convolution product and the norm $\|\cdot\|_{1,w}$ is a Banach algebra with $\Delta(L^{1}(G,w))=\widehat{G}(w)$, where $\widehat{G}(w)$ is the space of all non-zero complex-valued continuous homomorphisms $\varphi$ on $G$ such that $|\varphi(x)|\leq w(x)$ for each $x\in G$; see \cite{Kaniuth2}.

The space $M(G, w)$ of all complex regular Borel measures $\mu$ on $G$ such that $\mu w\in M(G)$  with convolution product and norm
 $$\|\mu\|_{M(G), w}=\|\mu w\|_{M(G)}=\int w(x)\ d|\mu|(x)$$
 is a Banach algebra called the weighted measure algebra, where $\mu w$ is defined by
$$\mu w(B)=\int_{B}w(x)\ d\mu(x)\quad \text{for each Borel subset $B$ of $G$  }.$$

%
\begin{proposition}\label{prop}
If $|\tau|\geq 1$ a.e., then
$L^{1}(G)_{\tau(n)}$ and $L^{1}(G, |\tau^{n}|)$ are topologically isomorphic.
\end{proposition}
\begin{proof}
Suppose that for almost every
$x\in G$, $|\tau(x)|\geq 1$.
Clearly if $f\in L^{1}(G)_{\tau(n)}$, then we have $f\in L^{1}(G, |\tau^{n}|)$. On the other hand, if $f\in L^{1}(G, |\tau^{n}|)$ by applying  $|\tau|\geq 1$ a.e., we conclude that $f\in L^{1}(G)_{\tau(n)}$, since for all $0\leq k\leq n$, $|f(x)\|\tau^{k}(x)|\leq |f(x)\|\tau^{n}(x)|$ a.e.. Therefore, $L^{1}(G, |\tau^{n}|)=L^{1}(G)_{\tau(n)}$ as two sets.
 Also,  using $|\tau|\geq 1$ a.e., we have
 $$\|f\|_{1, w}\leq \|f\|_{\tau(n)}\leq (n+1)\|f\|_{1, w},$$
 where $w=|\tau^{n}|$. So, two the norms $\|\cdot\|_{1,w}$ and $\|\cdot\|_{\tau(n)}$ are equivalent, which completes the proof.
\end{proof}
\begin{remark}\label{Rem1}
Note that $L^{1}(G)_{\tau(n)}$ and $L^{1}(G, |\tau^{n}|)$ are not equal in general. For example, let $G=R^{+}$ be the multiplicative group of all positive real numbers, $n=2$ and $\tau(x)=\frac{1}{x}$ for all $x\in R^{+}$. Clearly, $\tau$ is measurable and sub-multiplicative. Also, it is easily verified that $L^{1}(G)_{\tau(2)}\subseteq L^{1}(G, |\tau^{2}|)$. Now, take $0<\alpha <1$ and put $$f(x)=\left\{
                 \begin{array}{ll}
                   0, & \hbox{$0<x<1$;} \\
                   x^{\alpha}, & \hbox{$1\leq x$.}
                 \end{array}
               \right.
$$ One can easily check that $f$ is in $L^{1}(G, |\tau^{2}|)$ but it is not in $L^{1}(G)_{\tau(2)}$. So, $$L^{1}(G, |\tau^{2}|)\neq L^{1}(G)_{\tau(2)}.$$
Also, if we put $g(x)=\chi_{(0,1]}$, then $g\in L^{1}(G)$ but $g\notin L^{1}(G)_{\tau(2)}$. Hence, $$L^{1}(G)\neq L^{1}(G)_{\tau(2)}.$$
\end{remark}
\begin{remark}\label{Rem4.6}
  Although in general $L^{1}(G, |\tau^{k}|)\neq L^{1}(G)_{\tau(n)}$ for every integer $k$ with $0\leq k\leq n$, we do have
  $$\overline{L^{1}(G)_{\tau(n)}}^{\|\cdot\|_{1,|\tau^{k}|}}=L^{1}(G, |\tau^{k}|),\quad \quad (0\leq k\leq n).$$
  To see that this holds, observe that $C_{c}(G)$ is dense in $L^{1}(G, |\tau^{k}|)$ and  similar to \cite[Lemma 1.3.5 (i)]{Kaniuth2}, one can see that $C_{c}(G)\subseteq L^{1}(G)_{\tau(n)}$.
\end{remark}
\begin{remark}
  Let $K\subseteq G$ be a relatively compact neighborhood of $e$. Put
$$\mathcal{U}_{K}=\{U\subseteq K : \text{$U$ is a relatively compact neighborhood of $e$}\}.$$
For each $U\in \mathcal{U}_{K}$, let $f_{U}=\frac{\chi_{U}}{|U|}$, where $|U|$ denotes the Haar measure of $U$. Since $K$ is relatively compact by \cite[Lemma 1.3.3]{Kaniuth2}, there exists a positive real number $b$ such that $|\tau(x)|\leq b$ for all $x\in K$. So
$$\|f_{U}\|_{\tau(n)}\leq 1+b+\ldots +b^n\hspace{0.5cm}(U\in \mathcal{U}_{K}).$$
Also, similar to the group algebra case, for each $f\in L^{1}(G)_{\tau(n)}$ we have $\|f_{U}\ast f-f\|_{\tau(n)}\longrightarrow 0$ when $U$ tends to $\{e\}$. Therefore, $\{f_{U}\}_{U\in \mathcal{U}_{K}}$ is a BAI for $L^{1}(G)_{\tau(n)}$.
\end{remark}
\begin{remark}
Using the above remarks one can see that in general $L^{1}(G)_{\tau(n)}$ is not an abstract Segal algebra with respect to $L^{1}(G)$, because $L^{1}(G)_{\tau(n)}$ has a b.a.i and in general $L^{1}(G)\neq L^{1}(G)_{\tau(n)}$. But it is well-known that if $\mathcal{S}$ is an abstract Segal algebra with respect to $L^{1}(G)$ such that $\mathcal{S}$ has a BAI, then $\mathcal{S}=L^{1}(G)$.
\end{remark}

Suppose that $G$ is compact and $w$ is a weight on $G$. So, by Lemma 1.3.3 and Corollary 1.3.4 of \cite{Kaniuth2}, there exists positive real number $b$ such that $1\leq w(x)\leq b$ for all $x\in G$. Hence, using Proposition \ref{prop} and the proof of Theorem \ref{Th2}, we have the following corollary. Note that $"\cong"$ means topologically isomorphic.
\begin{corollary}
If $G$ is a compact group, then $L^{1}(G)_{\tau(n)}$ is a $\mathrm{BSE}$-algebra and we have the following relations:
$$L^{1}(G, |\tau^{n}|)\cong L^{1}(G)_{\tau(n)}\cong L^{1}(G).$$
\end{corollary}

\begin{remark}\label{Rem}
For every integer $k$ with $0\leq k\leq n$, we have,
$$L^{1}(G)_{\tau(n)}\subseteq L^{1}(G, |\tau^{k}|).$$
Therefore,
$$\widehat{G}\cup \widehat{G}(|\tau|)\cup\ldots \cup \widehat{G}(|\tau^{n}|)\subseteq \Delta(L^{1}(G)_{\tau(n)}).$$

\end{remark}
\begin{remark}
 Put $M(G)_{\tau(n)}=\bigcap_{k=1}^{n}M(G, |\tau^{k}|)$ and define the following norm:
 $$\|\mu\|_{\tau(n)}=\sum_{k=0}^{n}\left\|\mu\ |\tau^{k}|\ \right\|_{M(G)}\hspace{0.5cm}(\mu\in M(G)_{\tau(n)}).$$
A direct use of the convolution product shows that $(M(G)_{\tau(n)}, \|\cdot\|_{\tau(n)})$ is a normed algebra such that
\begin{equation}\label{1}
  M(G)_{\tau(n)}\subseteq \mathcal{M}(L^{1}(G)_{\tau(n)}).
\end{equation}

 We do not know whether the converse of the above inequality holds or not. Clearly, if $\tau$ is bounded, then $M(G)\cong M(G)_{\tau(n)}$ and so the converse of the above inequality holds by Wendel's Theorem.
\end{remark}
We end this paper with the following questions and conjecture.
\begin{description}
  \item[{\bf Conjecture}] We conjecture that Theorem \ref{Th3} is valid for every locally compact group.
  \item[{\bf Question}] Are the converse of Theorem \ref{Th2} and Relation \ref{1} hold?
\end{description}

{\bf Acknowledgments.} The authors would like to thank the referee for his/her suggestions and comments which improve the presentation of the paper especially, giving a shorter proof for Theorem \ref{seg1}.  The first named author of the paper supported partially by a grant from Gonbad Kavous University.

\bibliographystyle{amsplain}
\bibliography{References}

\end{document}